\documentclass{amsart}

\usepackage{amscd}

\newtheorem{thm}{Theorem}[section]
\newtheorem{prop}[thm]{Proposition}
\newtheorem{lem}[thm]{Lemma}

\theoremstyle{definition}
\newtheorem{dfn}[thm]{Definition}

\theoremstyle{remark}
\newtheorem{rem}{Remark}
\newtheorem*{acknowledgment}{Acknowledgments}

\newcommand{\C}{\mathbb{C}}
\newcommand{\N}{\mathbb{N}}
\newcommand{\R}{\mathbb{R}}
\newcommand{\Z}{\mathbb{Z}}

\newcommand{\bH}{\mathbb{H}}
\newcommand{\CB}{\mathbb{V}}

\newcommand{\G}{\mathcal{G}}
\renewcommand{\H}{\mathcal{H}}
\newcommand{\E}{\mathcal{E}}

\newcommand{\Ker}{\mathrm{Ker}}
\newcommand{\Hom}{\mathrm{Hom}}

\newcommand{\im}{\sqrt{\! - \! 1}}
\renewcommand{\d}{\partial}
\renewcommand{\l}{\ell}
\newcommand{\un}{\underline}

\def\til#1{ \tilde{#1} }
\def\SDC#1{ \Z({#1})^\infty_D }
\def\CSDC#1{ \Z({#1})^\infty_{D, \C}}

\def\abs#1{ \lvert{#1}\rvert }
\def\p#1{ \lVert{#1}\rVert }

\title{An analogue of the space of conformal blocks in $(4k+2)$-dimensions}

\author{Kiyonori Gomi
}

\address{
Graduate school of Mathematical Sciences, The University of Tokyo,
Komaba 3-8-1, Meguro-Ku, Tokyo, 153-8914 Japan.
}

\email{kgomi@ms.u-tokyo.ac.jp}

\date{}

\begin{document}

\begin{abstract}
Based on projective representations of smooth Deligne cohomology groups, we introduce an analogue of the space of conformal blocks to compact oriented $(4k+2)$-dimensional Riemannian manifolds with boundary. For the standard $(4k+2)$-dimensional disk, we compute the space concretely to prove that its dimension is finite. 
\end{abstract}

\maketitle



\section{Introduction}

As a fundamental ingredient, the \textit{space of conformal blocks} (or the space of vacua) in the Wess-Zumino-Witten model has been investigated by many physicists and mathematicians. While its construction usually appeals to representations of affine Lie algebras \cite{T-U-Y,U}, the formulation by means of representations of loop groups (\cite{Bry-M,S,W2}) provides schemes for generalizations.

\medskip

The theme of the present paper is an analogue of the space of conformal blocks in $(4k+2)$-dimensions. The idea of introducing such an analogue is to utilize \textit{smooth Deligne cohomology groups} (\cite{Bry,D-F,E-V}), or the groups of \textit{differential characters} (\cite{C-S}), instead of loop groups. In \cite{Go1,Go2}, some properties of smooth Deligne cohomology groups, such as projective representations, are studied. In a recent work of Freed, Moore and Segal \cite{F-M-S}, similar representations are also studied in a context of \textit{chiral (or self-dual) $2k$-forms} (\cite{W2}) on $(4k+2)$-dimensional spacetimes. 

\medskip

Our analogue of the space of conformal blocks is a vector space $\CB(W, \lambda)$ associated to a compact oriented $(4k+2)$-dimensional Riemannian manifold with boundary and an element $\lambda$ in a finite set $\Lambda(\d W)$. The finite set $\Lambda(\d W)$ is the set of equivalence classes of irreducible \textit{admissible representations} (\cite{Go2}) of the smooth Deligne cohomology group $\G(\d W) = H^{2k+1}(\d W, \SDC{2k+1})$. As will be detailed in the body of this paper (Section \ref{sec:analogue}), $\CB(W, \lambda)$ consists roughly of (dual) vectors in an irreducible representation realizing $\lambda$ which are invariant under actions of \textit{chiral (or self-dual) $2k$-forms}(\cite{F-M-S,W1}) on $W$.

\medskip

In the case of $k = 0$, we can interpret $\CB(W, \lambda)$ as the space of conformal blocks (or modular functor \cite{S}) based on representations of abelian loop groups. For example, we take $W$ to be the 2-dimensional disk $W = D^2$. In this case, $\G(S^1) = H^1(S^1, \SDC{1})$ is isomorphic to the loop group $LU(1)$. Irreducible admissible representations give rise to irreducible \textit{positive energy representations} (\cite{P-S}) of the loop group $LU(1)$ of level 2, which are classified by $\Lambda(S^1) \cong \Z_2$. Then the definition of $\CB(D^2, \lambda)$ can be read as:
$$
\CB(D^2, \lambda)
=
\{ 
\psi : \H_\lambda \to \C |\ 
\mbox{invariant under $\mathrm{Hol}(D^2, \C/\Z)$}
\},
$$
where $\H_\lambda$ is an irreducible representation corresponding to $\lambda$ on which the group $\mathrm{Hol}(D^2, \C/\Z)$ of holomorphic maps $f : D^2 \to \C/\Z$ acts densely and linearly through the ``Segal-Witten reciprocity law'' \cite{Bry-M,S,W2}.

\medskip

A property generally desired for $\CB(W, \lambda)$ is its finite-dimensionality. In the case of $k = 0$, there is a result of Segal regarding the property \cite{S}. The purpose of this paper is to prove that $\CB(W, \lambda)$ is finite-dimensional at least in the case where $W$ is the $(4k+2)$-dimensional disk $D^{4k+2} = \{ x \in \R^{4k+2} |\ \lvert x \rvert \le 1 \}$. Note that we have $\Lambda(S^{4k+1}) = \{ 0 \}$ for $k > 0$. Then we can show:

\begin{thm} \label{thm:main}
If $k > 0$, then $\CB(D^{4k+2}, 0) \cong \C$.
\end{thm}

The essential part of the proof is a fact about chiral $2k$-forms on $D^{4k+2}$, which we derive from \cite{I-T}. (See Section \ref{sec:disk} for detail.) The proof of Theorem \ref{thm:main} is applicable to the case of $k = 0$, and we have:
$$
\CB(D^2, \lambda) \cong
\left\{
\begin{array}{cc}
\C, & (\lambda = 0) \\
0. & (\lambda = 1)
\end{array}
\right.
$$
This result is consistent with the known fact about the dimension of the space of conformal blocks in the $U(1)$ Wess-Zumino-Witten model at level 2 (\cite{S,U}).

\medskip

The finite-dimensionality of $\CB(W, \lambda)$ for general $W$ remains open at present. A possible approach toward the issue is to generalize Segal's idea (p.431, \cite{S}), which should be examined in future studies.


\section{Analogue of the space of conformal blocks}
\label{sec:analogue}

In this section, we introduce the vector space $\CB(W, \lambda)$. For this aim, we summarize some results in \cite{Go1,Go2}. In particular, we review central extensions of smooth Deligne cohomology groups, a generalization of the Segal-Witten reciprocity law, and admissible representations.

\subsection{Central extension}

To begin with, we recall the definition of \textit{smooth Deligne cohomology} \cite{Bry,D-F,E-V}. For a non-negative integer $p$ and a smooth manifold $X$, the \textit{(complexified) smooth Deligne cohomology group} $H^*(X, \CSDC{p})$ is defined to be the hypercohomology of the following complex of sheaves on $X$:
$$
\CSDC{q} : \ 
\Z \longrightarrow
\un{A}^0_\C \overset{d}{\longrightarrow}
\un{A}^1_\C \overset{d}{\longrightarrow}
\cdots \overset{d}{\longrightarrow}
\un{A}^{q-1}_\C \longrightarrow
0 \longrightarrow \cdots,
$$
where $\Z$ is the constant sheaf, and $\un{A}^q_\C$ the sheaf of germs of $\C$-valued $q$-forms. We fix a non-negative integer $k$, and put $\G(X)_\C = H^{2k+1}(X, \CSDC{2k+1})$ for a smooth manifold $X$.

\begin{prop}[\cite{Go1}]
For a compact oriented $(4k+1)$-dimensional smooth manifold $M$ without boundary, there is a non-trivial central extension $\til{\G}(M)_\C$ of $\G(M)_\C$:
$$
\begin{CD}
1 @>>> \C^* @>>> \til{\G}(M)_\C @>>> \G(M)_\C @>>> 1.
\end{CD}
$$
\end{prop}

The central extension $\til{\G}(M)_\C$ is induced from the group 2-cocycle $S_{M, \C} : \G(M)_\C \times \G(M)_\C \to \C/\Z$ defined by $S_{M, \C}(f, g) = \int_M f \cup g$, where $\int_M$ and $\cup$ are the cup product and the integration in smooth Deligne cohomology.

\medskip

For a smooth manifold $X$, the smooth Deligne cohomology $H^1(X, \CSDC{1})$ is naturally isomorphic to $C^\infty(X, \C/\Z)$. Thus, if $k = 0$ and $M = S^1$, then we can identify $\G(S^1)_\C$ with the loop group $L\C^*$. In this case, $\til{\G}(S^1)_\C$ is isomorphic to $\widehat{L\C^*}/\Z_2$, where $\widehat{L\C^*}$ is the \textit{universal central extension} of $L\C^*$, (\cite{P-S}).

\subsection{A generalization of the Segal-Witten reciprocity law}

Let $W$ be a compact oriented $(4k+2)$-dimensional Riemannian manifold $W$ possibly with boundary. We denote by $A^{2k+1}(W, \C)$ the space of $\C$-valued $(2k+1)$-forms on $W$. The Hodge star operator $* : A^{2k+1}(W, \C) \to A^{2k+1}(W, \C)$ satisfies $** = -1$. Notice that, in general, the smooth Deligne cohomology $\G(X_\C) = H^{2k+1}(X, \CSDC{2k+1})$ fits into the following exact sequence. 
$$
0 \to H^{2k}(W, \C/\Z) \to \G(W)_\C \overset{\delta}{\to} 
A^{2k+1}(W, \C)_\Z \to 0,
$$
where $A^{2k+1}(W, \C)_\Z \subset A^{2k+1}(W, \C)$ is the subgroup consisting of closed integral forms. Using $*$ and $\delta$, we define the subgroups $\G(W)_\C^{\pm}$ in $\G(W)_\C$ by
$$ 
\G(W)^\pm_\C = 
\{ f \in \G(W)_\C |\ 
\delta(f) \mp \im * \delta(f) = 0\}.
$$
We call $\G(W)_\C^+$ the \textit{chiral subgroup}, since $2k$-forms $B \in A^{2k}(W, \C)$ such that $dB = i*dB$ are called \textit{chiral (or self-dual) $2k$-forms}. (See \cite{F-M-S,W1} for example.)

\begin{prop}[\cite{Go1}] \label{prop:reciprocity}
For a compact oriented $(4k+2)$-dimensional Riemannian manifold $W$ with boundary, the following map is a homomorphism:
$$
\til{r}^+ : \G(W)_\C^+ \longrightarrow \til{\G}(\d W)_\C, \quad
f \mapsto (f|_{\d W}, 1).
$$
\end{prop}

In the case of $k = 0$, $W$ is a Riemann surface. Since $\G(W)_\C^+$ is identified with the group of holomorphic functions $f : W \to \C/\Z$, Proposition \ref{prop:reciprocity} recovers the ``Segal-Witten reciprocity law''(\cite{Bry-M,S,W2}) for $\widehat{L\C^*}/\Z_2$.

\subsection{Admissible representations}

The group $\G(X)_\C = H^{2k+1}(X, \CSDC{2k+1})$ can be thought of as a complexification of the (real) smooth Deligne cohomology $\G(X) = H^{2k+1}(X, \SDC{2k+1})$ defined as the hypercohomology of the following complex of sheaves:
$$
\SDC{2k+1}: \ 
\Z \longrightarrow
\un{A}^0 \overset{d}{\longrightarrow}
\un{A}^1 \overset{d}{\longrightarrow}
\cdots \overset{d}{\longrightarrow}
\un{A}^{2k} \longrightarrow
0 \longrightarrow \cdots,
$$
where $\un{A}^q$ is the sheaf of germs of $\R$-valued $q$-forms. 

\medskip

For a compact oriented $(4k+1)$-dimensional Riemannian manifold $M$ without boundary, \textit{admissible representations} of $\G(M)$ are introduced in \cite{Go2}. An admissible representation $\rho : \G(M) \times \H \to \H$ of $\G(M)$ is a certain projective representation on a Hilbert space $\H$, and gives a linear representation $\til{\rho} : \til{\G}(M) \times \H \to \H$ of the central extension $\til{\G}(M)$ induced from the natural inclusion $\G(M) \subset \G(M)_\C$:
$$
\begin{CD}
1 @>>> U(1) @>>> \til{\G}(M) @>>> \G(M) @>>> 1 \\
@. @VVV @VVV @VVV @. \\
1 @>>> \C^* @>>> \til{\G}(M)_\C @>>> \G(M)_\C @>>> 1.
\end{CD}
$$
The set $\Lambda(M)$ of equivalence classes of irreducible admissible representations of $\G(M)$ is a finite set \cite{Go2}. For example, if $H^{2k+1}(M, \Z)$ is torsion free, then we can identify $\Lambda(M)$ with $H^{2k+1}(M, \Z_2)$. We write $(\til{\rho}_\lambda, \H_\lambda)$ for the linear representation of $\til{\G}(M)$ realizing $\lambda \in \Lambda$. 

\begin{prop}[\cite{Go2}] \label{prop:complex_extension}
Let $M$ be a compact oriented $(4k+1)$-dimensional Riemannian manifold $M$ without boundary. For $\lambda \in \Lambda(M)$, there exists an invariant dense subspace $\mathcal{E}_\lambda \subset \H_\lambda$, and the representation $\til{\rho}_\lambda : \til{\G}(M) \times \mathcal{E}_\lambda \to \mathcal{E}_\lambda$ extends to a linear representation $\til{\rho}_\lambda : \til{\G}(M)_\C \times \mathcal{E}_\lambda \to \mathcal{E}_\lambda$ of $\til{\G}(M)_\C$.
\end{prop}

We notice that $\til{\rho}_\lambda(f) : \E_\lambda \to \E_\lambda$ is generally unbounded, so that the action of $\til{\G}(M)_\C$ on $\E_\lambda$ does not extends to the whole of $\H_\lambda$.

\medskip

In the case of $k = 0$ and $M = S^1$, we can identify $\G(S^1)$ with the loop group $LU(1)$, which has $\G(S^1)_\C \cong L\C^*$ as a complexification. Admissible representations of $\G(S^1)$ give rise to positive energy representations of level 2. As is known \cite{P-S}, the equivalence classes of irreducible positive energy representations of $LU(1)$ of level 2 are in one to one correspondence with the elements in $\Lambda(S^1) \cong \Z_2$. A positive energy representation of $LU(1)$ extends to a representation of $L\C^*$ on an invariant dense subspace.

\subsection{Analogue of the space of conformal blocks}

We use Proposition \ref{prop:reciprocity} and Proposition \ref{prop:complex_extension} to formulate our analogue of the space of conformal blocks:

\begin{dfn}
Let $W$ be a compact oriented $(4k+2)$-dimensional Riemannian manifold with boundary. For $\lambda \in \Lambda(\d W)$, we define $\CB(W, \lambda)$ to be the vector space consisting of continuous linear maps $\psi : \E_\lambda \to \C$ invariant under the action of $\G(W)_\C^+$ through $\til{r}^+$:
\begin{equation*}
\begin{split}
\CB(W, \lambda) 
&= \Hom(\mathcal{E}_\lambda, \C)^{\mathrm{Im}\til{r}^+} \\
&= \{ 
\psi : \mathcal{E}_\lambda \to \C |\ 
\psi( \tilde{\rho}_\lambda(\til{r}^+(f)) v ) = \psi(v) \ 
\mbox{for $v \in \mathcal{E}_\lambda$ and $f \in \G(W)_\C^+$} 
\}.
\end{split}
\end{equation*}
\end{dfn}

Since the subgroup $\C^*$ in $\til{\G}(M)_\C = \G(M)_\C \times \C^*$ acts on $\mathcal{E}_\lambda$ by the scalar multiplication, we can formulate $\CB(W, \lambda)$ in terms of the projective representation $(\rho_\lambda, \mathcal{E}_\lambda)$ corresponding to $(\til{\rho}_\lambda, \mathcal{E}_\lambda)$:
\begin{align*}
\CB(W, \lambda) 
&= \Hom(\E_\lambda, \C)^{\mathrm{Im}r^+} \\
&= \{ 
\psi : \E_\lambda \to \C |\ 
\psi( \rho_\lambda(r^+(f)) v ) = \psi(v) \ 
\mbox{for $v \in \E_\lambda$ and $f \in \G(W)_\C^+$} 
\},
\end{align*}
where $r^+ : \G(W)_\C^+ \to \G(\d W)_\C$ is the restriction: $r^+(f) = f|_{\d W}$. 

\begin{rem}
One may wonder why we use representations of $\til{\G}(M)_\C$ on pre-Hilbert spaces to formulate $\CB(W, \lambda)$, instead of unitary representations of $\til{\G}(M)$ on Hilbert spaces. The reason is that we cannot introduce a counterpart of the chiral subgroup $\G(W)_\C^+$ to $\G(W)$. Notice, however, that we can formulate $\CB(W, \lambda)$ as
$$
\CB(W, \lambda) 
= \{ 
\psi : \H_\lambda \to \C |\ 
\psi( \rho_\lambda(r^+(f)) v ) = \psi(v) \ 
\mbox{for $v \in \mathcal{E}_\lambda$ and $f \in \G(W)_\C^+$}
\},
$$
because $\E_\lambda$ is dense in $\H_\lambda$.
\end{rem}


\section{Calculation of $\CB(D^{4k+2}, \lambda)$}
\label{sec:disk}

In this section, we prove Theorem \ref{thm:main}. As preparations for the proof, we review in some detail the construction of irreducible representations of Heisenberg groups in \cite{P-S}. We also study chiral $2k$-forms on $\R^{4k+2}$ by the help of results in \cite{I-T}.

\subsection{Representation of Heisenberg group}
\label{subsec:Heisenberg_representation}

For a compact oriented $(4k+1)$-dimensional Riemannian manifold $M$ without boundary, the group $\G(M)_\C$ admits the decomposition:
\begin{align*}
\G(M)_\C 
&\cong
(A^{2k}(M, \C)/A^{2k}(M, \C)_\Z) \times H^{2k+1}(M, \Z) \\
&\cong 
(\bH^{2k}(M, \C)/\bH^{2k}(M, \C)_\Z) \times
d^*(A^{2k+1}(M, \C)) \times H^{2k+1}(M, \Z),
\end{align*}
where $\bH^{2k}(M, \C)$ is the group of $\C$-valued harmonic $2k$-forms, $\bH^{2k}(M, \C)_\Z = \bH^{2k}(M, \C) \cap A^{2k}(M, \C)_\Z$ the subgroup of integral harmonic $2k$-forms, and $d^* : A^{2k+1}(M, \C) \to A^{2k}(M, \C)$ the formal adjoint of the exterior differential. Thus, in particular, if $M$ is such that $b^{2k+1}(M) = 0$, then $\G(M)_\C \cong d^*(A^{2k+1}(M, \Z))$. The representations $(\til{\rho}_\lambda, \E_\lambda)$ of $\til{\G}(M)_\C$ in Proposition \ref{prop:complex_extension} are built on a projective representation $(\rho, E)$ of $d^*(A^{2k+1}(M, \C)$. We review here the construction of $(\rho, E)$ following \cite{P-S}, and give a simple consequence.

\bigskip

As in \cite{Go2}, we define the Hermitian inner product $( \ , \ )_V$ on $d^*(A^{2k+1}(M, \C))$ by that induced from the Sobolev norm $\p{ \ \cdot \ }_s$ with $s = 1/2$. (Our convention is that $( \ , \ )_V$ is $\C$-linear in the first variable, which differs from that in \cite{P-S}.) On the completion $V_\C$ of $d^*(A^{2k+1}(M, \C))$, we define the linear map $J : V_\C \to V_\C$ by $J = \til{J}/\abs{\til{J}}$, where $\til{J} : d^*(A^{2k+1}(M, \C)) \to d^*(A^{2k+1}(M, \C))$ is the differential operator $\til{J} = *d$. Then $J$ is a complex structure compatible with $( \ , \ )_V$, and satisfies:
$$
(\alpha, J\overline{\beta})_V = \int_M \alpha \wedge d \beta
$$ 
for $\alpha, \beta \in d^*(A^{2k+1}(M, \C))$. By means of $J$, we decompose $V_\C$ into $V_\C = W \oplus \overline{W}$, where $J$ acts on $W$ and $\overline{W}$ by $\im$ and $-\im$, respectively.

Then we let $E = \C \langle \epsilon_\xi |\ \xi \in W \rangle$ be the vector space generated by the symbols $\epsilon_\xi$ corresponding to $\xi \in W$, and $\langle \ , \ \rangle : E \times E \to \C$ the Hermitian inner product $\langle \epsilon_\xi, \epsilon_\eta \rangle = e^{2(\xi, \eta)_V}$. For $v_+ \in W$ and $v_- \in \overline{W}$, we define $\rho(v_+ + v_-) : \ E \to E$ by
$$
\rho(v_+ + v_-) \epsilon_\xi = 
\exp\left(
- ( v_+, \overline{(v_-)} )_V
- 2 ( \xi, \overline{(v_-)} )_V \right)
\epsilon_{\xi + v_+}.
$$
We can verify $\rho(v)\rho(v') \epsilon_\xi = e^{\im (v, J \bar{v}')_V} \rho(v + v') \epsilon_\xi$ for $v, v' \in V_\C$, so that we have a projective representation $\rho : V_\C \times E \to E$. Because the group 2-cocycle $S_{M, \C}$ on $d^*(A^{2k+1}(M, \C))$ has the expression :
$$
S_{M, \C} (\alpha, \beta) = \int_M \alpha \wedge d\beta \mod \Z,
$$
we get the projective representation $\rho : d^*(A^{2k+1}(M, \C)) \times E \to E$. 

\medskip

In general, $\rho(\alpha) : E \to E$ is unbounded. However, if $\alpha$ belongs to the real vector space $d^*(A^{2k+1}(M))$ underlying $d^*(A^{2k}(M, \C))$, then $\rho(\alpha) : E \to E$ is isometric. Thus, $\rho(\alpha)$ extends to a unitary map on the completion $H = \overline{E}$ of $E$, and we have an irreducible projective unitary representation $\rho : d^*(A^{2k+1}(M)) \times H \to H$. As is shown in \cite{P-S}, we can identify $\overline{E}$ with a completion of the symmetric algebra $S(W)$ by the mapping $\epsilon_\xi \mapsto e^\xi = \sum_{j = 0}^\infty \xi^j/j!$.

\bigskip

The next lemma is a simple consequence from the above construction.

\begin{lem} \label{lem:Heisenberg_representation}
Let $(\rho, E)$ be as above.

(a) The vector space $\Hom(E, \C)^W$ is generated by the continuous linear map $\chi : E \to \C$ defined by $\chi(v) = \langle v, \epsilon_0 \rangle$:
$$
\Hom(E, \C)^W = \C \langle \chi \rangle.
$$

(b) For a dense subspace $U$ in $W$, we have $\Hom(E, \C)^U = \Hom(E, \C)^W$.
\end{lem}

\begin{proof}
To prove (a), we begin with proving the $W$-invariance of $\chi$. Notice that $\chi(\epsilon_\xi) = 1$ for all $\xi \in W$. For $f \in W$ and $v = \sum_j c_j \epsilon_{\xi_j} \in E$, we have:
\begin{align*}
\chi(v) &= \sum_j c_j \chi(\epsilon_{\xi_j}) = \sum_j c_j, \\
\chi((\rho(f)v) &= \sum_j c_j \chi( \rho(f)\epsilon_{\xi_j} ) 
= \sum_j c_j \chi(\epsilon_{\xi_j + f}) = \sum_j c_j.
\end{align*}
Hence $\chi$ is invariant under the action of $W$, and $\C \langle \chi \rangle \subset \Hom(E, \C)^W$. To see $\C \langle \chi \rangle \supset \Hom(E, \C)^W$, we show that any $\psi \in \Hom(E, \C)^W$ is of the form $\psi = c \chi$ for some $c \in \C$. For $v = \sum_j c_j \epsilon_{\xi_j} \in E$, the invariance of $\psi$ leads to:
\begin{equation*}
\begin{split}
\psi(v) 
&= \sum_j c_j \psi(\epsilon_{\xi_j}) 
= \sum_j c_j \psi(\rho(\xi_j) \epsilon_0)
= \sum_j c_j \psi(\epsilon_0) \\
&= \psi(\epsilon_0) \sum_j c_j 
= \psi(\epsilon_0) \chi(v).
\end{split}
\end{equation*}
If we put $c = \psi(\epsilon_0)$, then $\psi = c \chi$. For (b), it suffices to prove the inclusion $\Hom(E, \C)^U \subset \Hom(E, \C)^W$. So we will show $\psi \in \Hom(E, \C)^U$ is also invariant under $W$. For $f \in W$, there is a sequence $\{ f_n \}$ in $U$ converging to $f$. Notice that $\rho(\cdot)v : W \to E$ is continuous for $v \in E$. Now, we have:
$$
\psi(\rho(f)v) 
= \psi(\rho(\lim_{n \to \infty} f_n) v)
= \lim_{n \to \infty} \psi(\rho(f_n)v)
= \lim_{n \to \infty} \psi(v) = \psi(v),
$$
so that $\psi \in \Hom(E, \C)^W$.
\end{proof}

\begin{rem}
The key to Lemma \ref{lem:Heisenberg_representation} (b) is that the map $\rho(\cdot)v : W \to E$ is continuous for each $v \in W$. The representations $\til{\rho}_\lambda : \til{\G}(M)_\C \times \mathcal{E}_\lambda \to \mathcal{E}_\lambda$ in Proposition \ref{prop:complex_extension} have the same property \cite{Go2}.
\end{rem}

\subsection{Chiral $2k$-forms on $\R^{4k+2}$}

The Laplacian $\Delta = dd^* + d^*d$ preserves the subspace $d^*(A^{2k+1}(S^{4k+1}, \C))$. For an eigenvalue $\l$ of $\Delta$, we define $V_\l$ to be the following eigenspace:
$$
V_\l = 
\{
\beta \in d^*(A^{2k+1}(S^{4k+1}, \C)) |\ 
\Delta \beta = \l \beta
\}.
$$
The complex structure $J$, introduced in the previous subsection, preserves $V_\l$. (In particular, $J = *d/\sqrt{\l}$ on $V_\l$.) So we have the decomposition $V_\l = W_\l \oplus \overline{W}_\l$, where $J$ acts on $W_\l$ and $\overline{W}_\l$ by $\im$ and $-\im$, respectively. 

\begin{prop} \label{prop:chiral_form_plane}
There is the following relation of inclusion:
$$
\bigoplus_{\l} W_{\l}
\subset 
\mathrm{Im} \{ i^* : A^{2k}(\R^{4k+2}, \C)^+ \to A^{2k}(S^{4k+1}, \C) \}
\subset W,
$$
where $\bigoplus$ means the algebraic direct sum, $\l$ runs through all the distinct eigenvalues,  $i : S^{4k+1} \to \R^{4k+2}$ is the inclusion, and $A^{2k}(\R^{4k+2}, \C)^{\pm}$ are the spaces of chiral and anti-chiral $2k$-forms on $\R^{4k+2}$:
$$
A^{2k}(\R^{4k+2}, \C)^\pm =
\{
B \in A^{2k}(\R^{4k+2}, \C) |\ d B \mp \im *d B = 0
\}.
$$
\end{prop}

For the proof, we use some results shown by Ikeda and Taniguchi in \cite{I-T}. To explain the relevant results, we introduce some notations. Let $S^i(\R^{4k+1})$ and $\Lambda^p(\R^{4k+1})$ be the spaces of the symmetric tensors of degree $i$ and anti-symmetric tensors of degree $p$. We put $P^p_i = S^i(\R^{4k+1}) \otimes \Lambda^p(\R^{4k+1}) \otimes \C$, and regard $P^p_i$ as a subspace in $A^p(\R^{4k+1}, \C)$. We then define the vector spaces:
\begin{align*}
H^p_i     &= \mathrm{Ker} \Delta \cap \Ker d^* \cap P^p_i, \\
{}'\!H^p_i  &= \mathrm{Ker}d \cap H^p_i, \\
{}''\!H^p_i &= \mathrm{Ker} i\Big(r \frac{d}{d r} \Big) \cap H^p_i,
\end{align*}
where $i\left(r \frac{\d}{\d r} \right)$ is the contraction with the vector field $r \frac{d}{d r} = \sum_{j = 1}^{4k+2}x_j \frac{d}{dx_j}$. 

Notice that the standard action of $SO(4k+2)$ on $\R^{4k+2}$ makes ${}'\!H^p_i$ and ${}''\!H^p_i$ into $SO(4k+2)$-modules. Similarly, $V_\l$ is also an $SO(4k+2)$-module. From \cite{I-T} (Theorem 6.8, p.\ 537), we can derive:

\begin{prop}[\cite{I-T}] \label{prop:Ikeda_Taniguchi}
Let $\l_1 < \l_2 < \l_3 < \cdots$ be the sequence of distinct eigenvalues of $\Delta$ on $d^*(A^{2k+1}(S^{4k+1}, \C))$. For $i \in \N$, we have:

(a) The maps $i^* : A^{2k}(\R^{4k+2}, \C) \to A^{2k}(S^{4k+1}, \C)$ and $d : A^{2k}(\R^{4k+2}, \C) \to A^{2k+1}(\R^{4k+2}, \C)$ induce the following isomorphisms of $SO(4k+2)$-modules:
$$
\begin{CD}
V_{\l_i} 
@<{i^*}<<
{}''\!H^{2k}_i 
@>{d}>>
{}'\!H^{2k+1}_{i-1}.
\end{CD}
$$

(b) The $SO(4k+2)$-module $V_{\l_i}$ decomposes into two distinct irreducible modules having the same dimensions.
\end{prop}

\begin{rem}
More precisely, the sequence $\{ \l_i \}_{i \in \N}$ is given by $\l_i = (2k + i)^2$, and the dimension of the two irreducible modules in $V_{l_i}$ is $\binom{4k + i}{2k} \binom{2k + i -1}{2k}$. 
\end{rem}

We also note the next lemma for later use:

\begin{lem} \label{lem:chiral_or_achiral}
Let $( \ , \ )_{L^2}$ be the $L^2$-norm on $A^{2k+1}(D^{4k+2}, \C)$.

(a) For $B, B' \in A^{2k}(\R^{4k+2}, \C)$, we have:
$$
(i^*B, J i^*B')_V = - (dB, *dB')_{L^2}.
$$

(b) If $B \in A^{2k+1}(\R^{4k+2}, \C)$ obeys $(J - \sqrt{-1})i^*B = 0$, then:
$$
\p{H^+}^2_{L^2} - \p{H^-}^2_{L^2} \ge 0,
$$
where $H^{\pm} = (1 \pm \sqrt{-1}*)dB/2$. Similarly, if  $(J + \sqrt{-1})i^*B = 0$, then:
$$
\p{H^-}^2_{L^2} - \p{H^+}^2_{L^2} \ge 0.
$$
\end{lem}

\begin{proof}
We can readily show (a) combining properties of $( \ , \ )_V$ and $J$ with Stokes' theorem. Notice that the eigenspaces $\mathrm{Ker} (1 \pm \im *)$ in $A^{2k+1}(D^{4k+2}, \C)$ are orthogonal to each other with respect to the $L^2$-norm. Then the inequalities in (b) follow from $(i^*B, i^*B)_V \ge 0$ and (a).
\end{proof}

Proposition \ref{prop:Ikeda_Taniguchi} and the above lemma yield:

\begin{lem} \label{lem:chiral_extension_plane}
The map $i^*$ induces the following isomorphisms for $i \in \N$:
$$
{}''\!H^{2k}_i \cap A^{2k}(\R^{4k+2}, \C)^+ \cong W_{\l_i}, \quad
{}''\!H^{2k}_i \cap A^{2k}(\R^{4k+2}, \C)^- \cong \overline{W}_{\l_i}.
$$
\end{lem}

\begin{proof} 
Notice that the action of $SO(4k+2)$ on $V_{\l_i}$ is compatible with $J$. So $W_{\l_i}$ and $\overline{W}_{\l_i}$ are $SO(4k+2)$-modules. The dimensions of $W_{\l_i}$ and $\overline{W}_{\l_i}$ are the same, since they are complex-conjugate to each other. Similarly, since the $SO(4k+2)$-action on ${}'\!H^{2k+1}_{i-1}$ is compatible with the Hodge star operator $*$, the vector spaces $({}'\!H^{2k+1}_{i-1})^{\pm} = {}'\!H^{2k+1}_{i-1} \cap \mathrm{Ker}(1 \mp \im*)$ are also $SO(4k+2)$-modules with the same dimensions. Thus, by Proposition \ref{prop:Ikeda_Taniguchi}, $W_{\l_i}$ is isomorphic to one of $({}'\!H^{2k+1}_{i-1})^{\pm}$ through $d \circ (i^*)^{-1}$, and $\overline{W}_{\l_i}$ is isomorphic to the other. To settle the case, we appeal to Lemma \ref{lem:chiral_or_achiral} (b). Then the case of $W_{\l_i} \cong ({}'\!H^{2k+1}_{i-1})^+$ and $\overline{W}_{\l_i} \cong ({}'\!H^{2k+1}_{i-1})^-$ is consistent. Now the isomorphisms $d : {}''\!H^{2k}_i \cap A^{2k}(\R^{4k+2}, \C)^\pm \to ({}'\!H^{2k+1}_{i-1})^\pm$ complete the proof.
\end{proof}

\begin{proof}[The proof of Proposition \ref{prop:chiral_form_plane}]
By Lemma \ref{lem:chiral_extension_plane} we have:
$$
W_{\l_i} \subset 
\mathrm{Im} \{ i^* : A^{2k}(\R^{4k+2}, \C)^+ \to A^{2k}(S^{4k+1}, \C) \},
$$
which leads to the first inclusion in Proposition \ref{prop:chiral_form_plane}. For the second inclusion, we recall that the subspaces $W$ and $\overline{W}$ in $V_\C$ are orthogonal with respect to $( \ , \ )_V$. So, it suffices to verify the image $i^*(A^{2k}(\R^{4k+2}, \C)^+)$ is orthogonal to $\overline{W}$. By Lemma \ref{lem:chiral_extension_plane}, we also have:
$$
\overline{W}_{\l_i} \subset 
\mathrm{Im} \{ i^* : A^{2k}(\R^{4k+2}, \C)^- \to A^{2k}(S^{4k+1}, \C) \}.
$$
Thus, by the help of Lemma \ref{lem:chiral_or_achiral} (a), we see that $i^*(A^{2k}(\R^{4k+2}, \C)^+)$ is orthogonal to each $\overline{W}_{\l_i}$. Because $\bigoplus_i W_{\l_i}$ forms a dense subspace in $W$, the image $i^*(A^{2k}(\R^{4k+2}, \C)^+)$ is orthogonal to $\overline{W}$.
\end{proof}

\subsection{Proof of the main result}

We now compute $\CB(D^{4k+2}, \lambda)$.

First, we consider the case of $k > 0$. In this case, we have:
$$
\G(S^{4k+1})_\C 
= A^{2k}(S^{4k+1}, \C)/A^{2k}(S^{4k+1}, \C)_\Z 
\cong d^*(A^{2k+1}(S^{4k+1}, \C)).
$$
So the projective unitary representation $(\rho, H)$ reviewed in Subsection \ref{subsec:Heisenberg_representation} realizes the unique element in $\Lambda(S^{4k+1}) = \{ 0 \}$, and $E$ gives the invariant dense subspace in Proposition \ref{prop:complex_extension}.

\begin{thm}
If $k > 0$, then $\CB(D^{4k+2}, 0) \cong \C$.
\end{thm}

\begin{proof}
Note that $\G(D^{4k+2})^+_\C = A^{2k}(D^{4k+2}, \C)^+/A^{2k}(D^{4k+1}, \C)_\Z$.  Proposition \ref{prop:chiral_form_plane} leads to: $U \subset \mathrm{Im} r^+ \subset W$, where the dense subspace $U$ in $W$ is given by $U = \bigoplus_{i \in \N} W_{\l_i}$. This relation of inclusion implies:
$$
\Hom(E, \C)^U \supset \Hom(E, \C)^{\mathrm{Im} r^+} \supset \Hom(E, \C)^W.
$$
Therefore Lemma \ref{lem:Heisenberg_representation} establishes the theorem.
\end{proof}

In the case of $k = 0$, we have the familiar decomposition of $\G(S^1)_\C \cong L\C^*$:
$$
\G(S^1)_\C \cong \C/\Z \times 
\{ \phi : S^1 \to \R |\ \mbox{$\int \phi(\theta) d\theta = 0$} \} \times \Z.
$$
As is mentioned, admissible representations of $\G(S^1)$ are equivalent to positive energy representations of $LU(1)$ of level 2. For $\lambda \in \Lambda(S^1) = \Z_2 = \{ 0, 1 \}$, the invariant dense subspace $\E_\lambda$ in Proposition \ref{prop:complex_extension} is given by $\E_\lambda = \bigoplus_{\xi \in \Z}E_{\lambda + 2\xi}$, where $E_{\lambda + 2\xi} = E$ is the pre-Hilbert space in Subsection \ref{subsec:Heisenberg_representation} and the subgroup of constant loops $\C/\Z \subset \G(S^1)_\C$ acts on $E_{\lambda + 2\xi}$ with weight $\lambda + 2\xi$.

\begin{prop}
For $\lambda \in \Lambda(S^1) = \Z_2$, we have:
$$
\CB(D^2, \lambda)
\cong
\left\{
\begin{array}{cc}
\C, & (\lambda = 0) \\
0. & (\lambda = 1)
\end{array}
\right.
$$
\end{prop}

\begin{proof}
Clearly, constant loops $S^1 \to \C/\Z$ extend to holomorphic maps $D^2 \to \C/\Z$. So we use Proposition \ref{prop:chiral_form_plane} to obtain the relation of inclusion:
$$
\C/\Z \times U \subset \mathrm{Im} r^+ \subset \C/\Z \times W,
$$ 
where $U = \bigoplus_{i \in \N} W_{\l_i}$. Since $\C/\Z$ acts on $E_{\lambda + 2\xi}$ with weight $\lambda + 2\xi$, we have:
$$
\Hom(\E_\lambda, \C)^{\C/\Z} \subset
\prod_{\xi \in \Z} \Hom(E_{\lambda + 2\xi}, \C)^{\C/\Z} \cong
\left\{
\begin{array}{cc}
\Hom(E_0, \C), & (\lambda = 0) \\
\{ 0 \}. & (\lambda = 1)
\end{array}
\right.
$$
Thus, if $\lambda = 1$, then $\CB(D^2, \lambda) = \{ 0 \}$. In the case of $\lambda = 0$, we have:
$$
\Hom(E_0, \C)^U \supset 
\Hom(E_0, \C)^{\mathrm{Im}r^+} \supset 
\Hom(E_0, \C)^W.
$$
Now, Lemma \ref{lem:Heisenberg_representation} completes the proof.
\end{proof}


\bigskip

\begin{acknowledgment}
I would like to thank M. Furuta, T. Kohno and Y. Terashima for valuable discussions and helpful suggestions.
\end{acknowledgment}


\end{document}